\documentclass[a4paper,11pt]{amsart}

\usepackage[english]{babel} 
\usepackage{amsmath, enumerate, amsfonts,amssymb, amsthm}
\usepackage[dvips]{graphicx}

\title{Local tropical variety}
\author{Naoyuki Touda}
\address{Department of Mathematics, Kobe University,
Rokko 1-1, Kobe 657-8501, Japan}
\email{touda@math.kobe-u.ac.jp }

\newtheorem{definition}{\bf Def\mbox{}inition}[section]
\newtheorem{theorem}[definition]{\bf Theorem}
\newtheorem{proposition}[definition]{\bf Proposition}
\newtheorem{lemma}[definition]{\bf Lemma}

\newtheorem{example}[definition]{\bf Example}

\newcommand{\C}{\mathbb{C}}
\newcommand{\N}{\mathbb{N}}

\newcommand{\gr}{\mathrm{Gr}}   
\newcommand{\ini}{\mathrm{in}}  
\newcommand{\lc}{\mathrm{lc}}   

\newcommand{\Ohat}{{\hat{\mathcal{O}}}}

\newcommand{\Uloc}{\mathcal{U}_{\mathrm{loc}}} 

\begin{document}

\maketitle

\section{Introduction}
The tropical variety of an ideal of $K[x] := K[x_1,\ldots,x_n]$, where $K = \mathbb{C} \{ \{ t \} \} $ is 
the field of Puiseux series, was introduced by D.Speyer and B.Sturmfels \cite{Spey}. 
Let us sketch their construction. Any Puiseux series $p(t)$ can be written as  $p(t) = c_1t^{q_1}+c_2t^{q_2}+c_3t^{q_3}+\cdots$, 
where $c_1, c_2, \cdots$ are non-zero complex numbers and $q_1 < q_2 < \cdots$ are rational numbers 
with common denominator.
The order of $p(t)$, denoted by $\mbox{ord}(p(t))$, is the exponent $q_1$. 
Then we have the following map:
$$\begin{array}{cccc} 
     \mbox{ord}   : & (K \setminus \{0\})^n  & \to   & \mathbb{Q}^n   \subset \mathbb{R}^n \\
			   & \mbox{ord}(p_1(t),\ldots,p_n(t)) & \mapsto  &  (\mbox{ord}((p_1(t)), \ldots,\mbox{ord}(p_n(t))))
\end{array} $$
We fix a weight vector $w \in \mathbb{R}^n$. The weight of the variable $x_i$ is $w_i$. 
The weight of a term $p_i(t)x_1^{a_{i1}}\cdots x_n^{a_{in}}$ is the real number 
$\mbox{\rm ord}(p_i(t))+a_{i1}w_1+\dots a_{in}w_n$.
Consider a polynomial $f = \sum p_i(t)x_1^{a_{i1}}\cdots x_n^{a_{in}} \in K[x]$. 
Let $\overline{w}$ be the smallest weight among all the terms in $f$. The initial form of $f$ equals 
$$\mbox{\rm in}_w(f) \ = \ \sum  c_{\alpha_1,\ldots,\alpha_n} \cdot x_1^{\alpha _1}\cdots x_n^{\alpha _n}$$
where the sum ranges over all the terms $p_i(t)x_1^{a_{i1}}\cdots x_n^{a_{in}}$ in $f$
whose $w$-weight coincides with $\overline{w}$ and where $c_{\alpha_1,\ldots,\alpha_n} \in \mathbb{C}$ denotes 
the coefficient of $t^{\overline{w} - \alpha_1w_1 - \cdots - \alpha_nw_n }$
in the Puiseux series $p(t)$. 
The initial ideal $\mbox{\rm in}_w(I) \subset K[x]$ is defined as the ideal 
generated by all the initial forms $\mbox{\rm in}_w(f)$ where $f$ runs over $I$.

Given an ideal $I$ in $K[x_1,\ldots,x_n]$. We define its zero set 
$V(I):=  \{ v \in  (K \setminus \{0\})^n  ~|~  f(v)=0 , {}^\forall f \in I \}$. 
Then the tropical variety $\mathcal{T}(I)$ is defined as 
$\mathcal{T}(I) \ = \ \overline{ \mbox{\rm ord}(V(I))} \subset \mathbb{R}^n$, the topological closure of the image of $V(I)$ under the map above.

In tropical geometry we consider the tropical semiring 
$(\mathbb{R} \cup \{ +\infty \} , \oplus, \odot)$, which has tropical addition $\oplus$ and 
tropical multiplication $\odot$, that means, for $x,y \in \mathbb{R} \cup \{ +\infty \}, \ 
x \ \oplus \ y \ := \ \mbox{\rm min} \{ x,y \} \mbox{ and } x \ \odot \ y \ := \ x + y$. 
Then for vectors $(a_1,\ldots,a_n), ~(b_1,\ldots,b_n) \in \mathbb{R}^n$, its tropical addition is 
$(a_1,\ldots,a_n) \oplus (b_1,\ldots,b_n) := (\mbox{\rm min} \{ a_1,b_1 \}, \ldots, \mbox{\rm min} \{ a_n,b_n \} )$
and its tropical scalar multiplication is 
$\lambda \odot (a_1,a_2,\ldots,a_n) := (\lambda + a_1,\lambda + a_2,\ldots,\lambda + a_n)$ ($\lambda \in \mathbb{R}$).
A tropical monomial is an expression of the form 
$c\odot x_1^{a_1} \odot \cdots  \odot  x_n^{a_n}, ~c\in \mathbb{R}$ 
where the powers of the variables are computed tropically as well, 
for instance, $x_1^3 = x_1 \odot x_1 \odot x_1$.
A tropical polynomial is a finite tropical sum of tropical monomials, 
$F  =   c_1\odot x_1^{a_{11}} \odot \cdots  \odot  x_n^{a_{1n}} \oplus \cdots \oplus c_r \odot x_1^{a_{r1}} \odot \cdots  \odot  x_n^{a_{rn}}$.
Then for any $(w_1,\ldots,w_n) \in \mathbb{R}^n$,
$F(w_1,\ldots,w_n) =  \mbox{\rm min}_{1 \leq i \leq r}  \{ c_i+a_{i1}w_1+\cdots +a_{in}w_n \}$. 
The tropical hypersurface of $F$, denoted by  $\mathcal{T}_h(F)$,  is the set of  points 
$(w_1,\ldots,w_n) \in \mathbb{R}^n$ at which the minimum of $F(w_1,\ldots,w_n)$ is attained twice or more.
For a polynomial $f =\sum_{i=1}^r p_i(t)x_1^{a_{i1}}\dots x_n^{a_{in}} \ (p_i(t) \in K)$ in $K[x]$, 
$\mbox{\rm trop }(f) \ := \ c_1\odot x_1^{a_{11}} \odot \cdots  \odot  x_n^{a_{1n}} \oplus \cdots \oplus c_r \odot x_1^{a_{r1}} \odot \cdots  \odot  x_n^{a_{rn}}, \ (c_i = \mbox{ \rm ord}(p_i(t)))$
is called the tropicalization of $f$.

\

With the notations above, the following theorem was proved by D.Speyer and B.Sturmfels:
\begin{theorem}\label{main}\rm{[\cite{Spey} Theorem 2.1.]}  For any ideal $I$ in $ K[x]$, the following subsets of
$\mathbb{R}^n$  coincide:
\begin{enumerate}
\item $\mathcal{T}(I)=\overline{\mbox{\rm ord}(V(I))};$

\item $\bigcap_{f\in I} \mathcal{T}_h(\mbox{trop}(f));$

\item $\{ w\in \mathbb{R}^n \ | \ \mbox{for any }  f \in I, \ \mbox{\rm in}_w(f) \mbox{ is not a monomial} \}.$
\end{enumerate}
\end{theorem}

In the particular case $K=\mathbb{C}$, we have:

\begin{proposition}\label{subfan}
Let $I$ be a homogeneous ideal in $\mathbb{C}[x]$, then $\mathcal{T}(I)$ is a subfan of the 
Gr\"obner fan $\mbox{\rm GF}(I)$.
\end{proposition}

In this paper, we will give an analogous definition of local tropical varieties  and prove analogous theorems  in  the formal power series ring 
$\mathbb{C}[[x_1,\ldots,x_n]] = \hat{\mathcal{O}}$. 

\section{Local tropical hypersurface}

In this section we will define local tropical hypersurface. 
First, let's introduce some notations which are necessary for our local construction.

We denote by $\mathcal{U}_{\mbox{\rm loc}}$ the set of local weights
$\{ u ~\in ~ \mathbb{R}^n~|~u_i \geq 0, \forall i \}$.
For $f = \sum_{\alpha}p_{\alpha}x^{\alpha} \in \Ohat$, we define it's support 
$\mbox{\rm Supp}(f)\subset \mathbb{N}^n$ 
as the set of $\alpha$ such that $p_{\alpha}\neq 0$.
For a local weight vector $u \in \mathcal{U}_{\mbox{\rm loc}}$, the $u$-weight of $f$ (denoted by  
$\mbox{\rm wei}^u(f)$) is the minimum of the scalar products 
$u \cdot \alpha $ for $u \in \Uloc$. The weight gives rise to a filtration 
$F^u(\Ohat)$ given by $F_k^u : = \{ f \ | \ \mbox{\rm wei}^u(f) \geq k \}$ as well as to the associated 
graded ring $\mbox{\rm Gr}^u(\Ohat):= \bigoplus_{k} F_k^u / F^u_{k > }$ (direct sum).
For $f \in \Ohat$, its initial form $\mbox{\rm in}_u(f)$ shall be the class of $f$ in 
$F_k^u / F^u_{k > }$ where $k = \mbox{\rm wei}^u(f)$.
If $I \subset \Ohat$ is an ideal, then its initial ideal $\mbox{\rm in}_u(I)$ is an ideal in 
$\mbox{\rm Gr}^u(\Ohat)$
generated by all the $\ini_u(f)$ for $f\in I$, in other words,  
$\mbox{\rm in}_u(I):= \langle \mbox{\rm in}_u(f) \ | \ f \in I \rangle_{\mbox{\rm Gr}^u(\Ohat)}$.

\begin{lemma}\label{3-13-1}
Let $u \in \Uloc$ be a local weight vector and suppose that
$u=(u_1,\ldots,u_m, 0, \ldots ,0)$
with $0\le m\le n$ and $u_i > 0$. Then 
$$\gr^u(\Ohat)=\mathbb{C} [x_1, \ldots,x_m][[ x_{m+1}, \ldots, x_n]].$$
\end{lemma}

Thanks to the preceding lemma, we shall see $\gr^u(\Ohat)$ as a subring of $\Ohat$.

\begin{definition}  \rm
For $f = \sum_{\alpha}p_{\alpha}x^{\alpha} \in \Ohat$, 
 we define the tropicalization of $f$ as 
$$\mbox{\rm trop}(f) := \bigoplus_{\alpha \in \mbox{\rm Supp}(f)}a_{\alpha} \odot x^{\alpha} \mbox{ \ (infinite tropical sum)}$$
and its tropical hypersurface as
$$\mathcal{T}^h_{\mbox{\rm loc}}(\mbox{\rm trop}(f)) :=\{ w \in \mathcal{U}_{\mbox{\rm loc}} \ | \ \mbox{the minimum is attained twice or more} \}. $$

\end{definition}

The following proposition, which is analogous to the case of the polynomial ring, holds:

\begin{proposition}\label{2-1-5} For any ideal $I \subset \Ohat$, the following subsets of $\Uloc$ coincide:
\begin{enumerate}
\item $\{ w \in \mathcal{U}_{\mbox{\rm loc}} \ | \ \forall f \in I, \ \mbox{\rm in}_w(f) \in  \mbox{\rm Gr}^w(\Ohat)  \mbox{ is not a monomial}  \};$ 

\item $\bigcap_{f \in I} \mathcal{T}^h_{\mbox{\rm loc}}(\mbox{\rm trop}(f)).  $
\end{enumerate}

\end{proposition}

\begin{proof}
The proof is analogous to the case of polynomial ring,  \cite[Theorem 2.1]{Spey} .
\end{proof}

\section{Local tropical variety (principal case)}

In this section we will define the local tropical variety of a principal ideal in $\Ohat$. For this purpose,
the following proposition will play an important role.

\begin{proposition}\label{E(f)}  For $f = \sum_{\alpha \in \mbox{\rm Supp}(f)}p_{\alpha}x^{\alpha} \in \mathbb{C}[[x]]$, 
there exists the finite and minimal subset 
$E_0(f) \subset \mbox{ \rm Supp}(f) $ s.t. 
$$E_0(f) + \mathbb{N}^n = \mbox{\rm Supp}(f) + \mathbb{N}^n. $$
\end{proposition}

Now we will construct the local tropical variety of a principal ideal in $\widehat{\mathcal{O}}$.

\subsection{Tropical Variety on the Maximal Stratum}

In this subsection, we suppose $w$ lies in  $\mathcal{U}_{\mbox{ \rm loc}}^0 :=\{ (u_1,\ldots,u_n) \ | \ {}^\forall i, u_i > 0 \} $.

\begin{definition}\label{25-0}  
Let $f = \sum_{\alpha \in \mbox{\rm Supp}(f)}p_{\alpha}x^{\alpha} \in \mathbb{C}[[x]]$, then we define 
$$\widetilde{f}^0 := \sum_{\beta \in E_0(f) }p_{\beta}x^{\beta} \in \mathbb{C}[x] ( \subset \mathbb{C}[[x]]).$$
\end{definition}

\begin{proposition}\label{25-1} The following subsets of $\mathcal{U}_{\mbox{loc}}^0$  coincide:
\begin{enumerate}
\item $\{ w \in \mathcal{U}_{\mbox{loc}}^0 \ | \ \forall f \in I, \ \mbox{\rm in}_w(f) \in  \mbox{\rm Gr}^w(K[[x]])  \mbox{ is not a monomial}   \};$

\item $\mathcal{T}^h_{\mbox{\rm loc}}(\mbox{\rm trop}(f)) \bigcap \mathcal{U}_{\mbox{\rm loc}}^0 ;$

\item $\mathcal{T}^h_{\mbox{\rm loc}}(\mbox{\rm trop}(\widetilde{f}^0)) \bigcap \mathcal{U}_{\mbox{\rm loc}}^0 ;$
\item $\overline{\mbox{\rm ord}(V( \langle \widetilde{f}^0 \rangle ) )} \bigcap \mathcal{U}_{\mbox{\rm loc}}^0 .$
\end{enumerate}
Here, $V( \langle \widetilde{f}^0 \rangle )$
is defined by 
$\{ v \in (\mathbb{C} \{ \{ t \} \} \setminus \{ 0 \})^n  \ | \ \widetilde{f}^0(v)=0 \}$.
\end{proposition}

To prove the previous proposition, we will use the next lemma:

\begin{lemma}\label{25-33}  
Fix $u=(u_1,\ldots,u_n) \in \mathcal{U}_{\mbox{\rm loc}}^0$. 
Suppose that $(\alpha_1,\ldots,\alpha_n) \in \mbox{\rm Supp}(f)$ satisfies 
$\alpha_1u_1+ \cdots + \alpha_n u_n \leq 
 \alpha_1'u_1 + \cdots + \alpha_n'u_n$ 
for any 
$ ( \alpha_1', \ldots, \alpha_n') \in \mbox{\rm Supp}(f)$.
Then $(\alpha_1,\ldots,\alpha_n)$ is in $ E_0(f)$. 
\end{lemma}
The proof is easy.

\begin{proof}[Proof of Proposition \ref{25-1}]

Proofs of (1)$=$(2) and (3)$=$(4) is similar to the case of polynomial ring. See \cite[Theorem 2.1.]{Spey}.

((2) $\subset$ (3)) :
Take $w=(w_1,\ldots,w_n) \in$ (2).
Then there exist $\alpha_1,\alpha_2 \in \mbox{Supp}(f) \mbox{ satisfying } \alpha_1 \cdot w = \alpha_2 \cdot w \leq \alpha \cdot w \mbox{ for any } \alpha \in \mbox{Supp}(f)$.
By \mbox{Lemma } \ref{25-33}, this implies $\alpha_1,\alpha_2 \in E_0(f)$. 
 Since $ E_0(f) \subset \mbox{Supp}(f)$, $w \in$ (3)

((3) $\subset$ (2)) :
Take $w=(w_1,\ldots,w_n) \in$ (3). 
Then there exist $\beta_1,\beta_2 \in E_0(f) \mbox{ satisfying } \beta_1 \cdot w = \beta_2 \cdot w \leq \beta \cdot w \mbox{ for any } \beta \in E_0(f)$.
By Lemma \ref{25-33}, this implies $\beta_1 \cdot w = \beta_2 \cdot w \leq \alpha \cdot w \mbox{ for any } \alpha \in \mbox{Supp}(f)$.
\end{proof}

\subsection{Tropical Variety on a substratum $\mathcal{U}_{\mbox{\rm loc}}^1$}

In this subsection, we suppose that $w$ lies in $ \mathcal{U}_{\mbox{ \rm loc}}^1 :=\{ (0,u_2,\ldots,u_n) \ | \ \forall i(\neq 1), u_i > 0 \}$.

First, we consider the following map:
$$\begin{array}{cccc} 
       \phi_1 : & \Ohat = \mathbb{C}[[x_1,\ldots,x_n]] & \to  & \bigl(\mathbb{C}[[x_1]]\bigr)[[x_2,\ldots,x_n]]. \\
			    &     f                & \mapsto & \phi_1(f)
\end{array} $$
That means, we consider $x_1$ as a parameter and $f=\phi_1(f)$ as an element of 
$\bigl(\mathbb{C}[[x_1]]\bigr)[[x_2,\ldots,x_n]]$, 
the formal power series ring with the variables $x_2,\ldots,x_n$ over the formal power series ring $\mathbb{C}[[x_1]]$.
Then we define $E_1(\phi_1(f)) \subset \mathbb{N}^{n-1}$ as in Proposition \ref{E(f)}.
We denote by $\Phi_1$ the projection map 
from the set of points $\mbox{Supp}(f)$ to
$\mathbb{N}^{n-1} \bigcap \mbox{Supp}(f)$.
Define  $E_1(f)  ( \subset \mbox{Supp}(f) )$
with $\Phi_1^{-1}(E_1(\phi_1(f)))$.

\begin{definition}\label{25--1} 
$\widetilde{f}^1 := \sum_{\beta \in E_1(f) }a_{\beta}x^{\beta} \in \Ohat$.
\end{definition}

Then, as an analogy of Lemma \ref{25-33}, we have:
\begin{lemma}\label{25-3}  Fix $u=(0,u_2,\ldots,u_n) \in \mathcal{U}^1_{\mbox{\rm loc}}$. 
Suppose that $(\alpha_1,\ldots,\alpha_n) \in \mbox{\rm Supp}(f)$ satisfies 
$\alpha_1 \cdot 0 + \alpha_2\cdot u_2+ \cdots + \alpha_n \cdot u_n \leq 
\alpha_1'\cdot 0 + \alpha \cdot u_2 + \cdots + \alpha_n' \cdot u_n$ 
 for any $ ( \alpha_1', \ldots, \alpha_n') \in \mbox{\rm Supp}(f)$.  
Then 
$(\alpha_1,\ldots,\alpha_n)$ is in $ E_1(f)$
\end{lemma}

Now, we suppose $E_1(\phi(f)) = \{ \beta^1=(\beta^1_2,\ldots,\beta^1_n), \ldots,  \beta^l=(\beta^l_2,\ldots,\beta^l_n) \}$.
Then $ E_1(f)  = \{ E_1^1,\ldots,E_1^l \} $, where $\Phi_1^{-1}(\beta^i)  =E_1^i$ ($1 \leq i \leq l$).

\begin{lemma}\label{25-2}

Take $\alpha^1,\alpha^2 \in E_1^i, (i=1,\ldots,l)$. For any $w \in  \mathcal{U}_{\mbox{\rm loc}}^1$,
 $\alpha^1 \cdot w = \alpha^2 \cdot w$.
\end{lemma}

Now, we set $ E_1(f) = \{E_1^1,\ldots, E_1^l \} $ as above. 
For each $E^i$, if $\sharp( E_1^i) \geq 2$, we choose two elements 
$\alpha^1,\alpha^2 \ (\alpha^1 \neq \alpha^2)$ from $E_1^i$ and we set $A_1^i = \{ \alpha^1,\alpha^2 \}$. 
If $\sharp( E_1^i)=1$, we set $A_1^i = E^i_1$. Then, we define $A_1(f)= \{ A_1^1,\ldots,A_1^l \} ( \subset E_1(f) )$.
Clearly, $A_1(f)$ is a finite set ($\sharp(A_i(f)) \leq 2l$).
 
\begin{definition}\label{2-2-16}

$\widehat{f}^1 : = \sum_{\alpha \in A_1(f)} a_{\alpha} x^{\alpha} \in \mathbb{C}[x]$ 

\end{definition}

\begin{proposition}\label{2-3-2}

$\mathcal{T}^h_{\mbox{\rm loc}}(\mbox{\rm trop}(\widetilde{f}^1)) \bigcap \mathcal{U}_{\mbox{\rm loc}}^1 = \mathcal{T}^h_{\mbox{\rm loc}}(\mbox{\rm trop}(\widehat{f}^1)) \bigcap \mathcal{U}_{\mbox{\rm loc}}^1  $
\end{proposition}

\begin{proof}

($\subset$): Take $w \in $ LHS. Then there exist  $b \in E_1^i \mbox{ and } c \in E_1^j (1 \leq i,j \leq l)$
 s.t. $b\cdot w = c \cdot w \leq a \cdot w \mbox{ for any } a \in E_1(f) $.
By Lemma \ref{25-2}, there exist $\hat{b} \in A_1^i \mbox{ and } \hat{c} \in A_1^j $ 
s.t. $b\cdot w = \hat{b}\cdot w, ~ c\cdot w = \hat{c}\cdot w$. 
Since $A_1(f) \subset E_1(f)$, this implies $w \in $ RHS . 

($\supset$): Take $w \in $ RHS. Then there exist $b \in A_1^i \mbox{ and } c \in  A_1^j (1 \leq i,j \leq l)$
 s.t. $b\cdot w = c \cdot w \leq a \cdot w \mbox{ for any } a \in A_1(f) $.
 Now we suppose that there exists $d \in E_1^k (\subset E_1(f)) \mbox{ s.t. } d \cdot w \lneq b\cdot w = c \cdot w$. 
Then, by Lemma \ref{25-2}, for any $\hat{d} \in A_1^k, \hat{d} \cdot w  = d \cdot w \lneq b \cdot w = c \cdot w$. 
This contradicts the hypothesis. 
Thus, $b\cdot w = c \cdot w \leq a \cdot w \mbox{ for any } a \in  E_1(f) $ and this implies $w \in $ LHS .
\end{proof}

Analogously to the Proposition \ref{25-1}, the following proposition holds:

\begin{proposition}\label{2-10-7} The following subsets of $\mathcal{U}_{\mbox{loc}}^1 $ coincide:
\begin{enumerate}
\item $\{ w \in \mathcal{U}_{\mbox{\rm loc}}^1 \ | \ \forall f \in I, \ \mbox{\rm in}_w(f) \in  \mbox{\rm Gr}^w(\Ohat)  \mbox{ is not a monomial} \};$ 

\item $\mathcal{T}^h_{\mbox{\rm loc}}(\mbox{\rm trop}(f)) \bigcap \mathcal{U}_{\mbox{\rm loc}}^1 ;$

\item $\mathcal{T}^h_{\mbox{\rm loc}}(\mbox{\rm trop}(\widetilde{f}^1)) \bigcap \mathcal{U}_{\mbox{\rm loc}}^1 ;$

\item $\mathcal{T}^h_{\mbox{\rm loc}}(\mbox{\rm trop}(\widehat{f}^1)) \bigcap \mathcal{U}_{\mbox{\rm loc}}^1 ;$

\item $\overline{\mbox{\rm ord}(V( \langle \widehat{f}^1 \rangle ) )} \bigcap \mathcal{U}_{\mbox{\rm loc}}^1 .$
\end{enumerate}
\end{proposition}

\begin{proof}
The proofs of (1)=(2) and (4)=(5) are analogous 
to the proof of the polynomial case. 
By Proposition \ref{2-3-2}, 
it is enough to prove the equality of (2) and (3).

((2) $\subset$ (3)): 
Take $w=(0,w_2,\ldots,w_n) \in$ (2). 
Then there exist
$\alpha_1,\alpha_2 \in \mbox{Supp}(f) \mbox{ s.t. } \alpha_1 \cdot w = \alpha_2 \cdot w \leq \alpha \cdot w \mbox{ for any } \alpha \in \mbox{Supp}(f)$.
By Lemma \ref{25-2}, this implies
$\alpha_1,\alpha_2 \in E_1(f)$. 
Since $ E_1(f) \subset \mbox{Supp}(f)$,
we have $ w \in$ (3).

((3) $\subset$ (2)) :
Take $w=(0,w_2,\ldots,w_n) \in (3) $
Then there exist 
$\beta_1,\beta_2 \in E_1(f) \mbox{ s.t. } \beta_1 \cdot w = \beta_2 \cdot w \leq \beta \cdot w \mbox{ for any } \beta \in E_1(f)$.
By Lemma \ref{25-2}, $ \beta_1 \cdot w = \beta_2 \cdot w \leq \alpha \cdot w \mbox{ for any } \alpha \in \mbox{Supp}(f) $.
This implies $w \in $ (2).
\end{proof}

\subsection{Tropical Variety on General Strata}

In this subsection, we suppose that $w$ lies in $\mathcal{U}_{\mbox{\rm loc}}$.

We will use the same argument as in 3.2. For 
$\mathcal{U}_{\mbox{loc}}^j :=\{ (u_1,\ldots,u_n) \ | \ \forall i ( \neq j ), u_i > 0, u_j=0\}, \ldots, \mathcal{U}_{\mbox{loc}}^n, \mathcal{U}_{\mbox{loc}}^{12} :=\{ (0,0,u_3\ldots,u_n) \ | \ \forall i ( \neq 1,2 ), u_i > 0 \},\ldots,\mathcal{U}_{\mbox{loc}}^{234\cdots n} :=\{ (0,\ldots,0,u_n) \ | \  u_n > 0, \forall j ( \neq n ),u_j=0\}, \mathcal{U}_{\mbox{loc}}^{123\cdots n} :=\{ (0,\ldots,0) \}$ .

\begin{definition}  For $\mathcal{U}_{\mbox{loc}}^{123\cdots n} :=\{ (0,\ldots,0) \}$,
we define $\widetilde{f}^{12\cdots n}=\widehat{f}^{12\cdots n} : = 0 $.
For $\mathcal{U}_{\mbox{loc}}^{0} :=\{ (u_1,\ldots,u_n) \ | \ \forall i, u_i > 0 \}$, 
we define $\widetilde{f}^{0}=\widehat{f}^{0} $.
\end{definition}
Then, for each $\lambda \in \Lambda = \{ 0,1,2,\ldots,n,12,13,\ldots,234\cdots n,123\cdots n \}$, 
the following proposition holds:

\begin{proposition}\label{2-2-1} The following subsets of  $\mathcal{U}_{\mbox{\rm loc}}^{\lambda}$ coincide:

\begin{enumerate}
\item $\{ w \in \mathcal{U}_{\mbox{\rm loc}}^{\lambda} \ | \ \forall f \in I, \ \mbox{\rm in}_w(f) \in  \mbox{\rm Gr}^w(\Ohat)  \mbox{ is not a monomial} \};$

\item $\mathcal{T}^h_{\mbox{\rm loc}}(\mbox{\rm trop}(f)) \bigcap \mathcal{U}_{\mbox{\rm loc}}^{\lambda} ;$

\item $\mathcal{T}^h_{\mbox{\rm loc}}(\mbox{\rm trop}(\widetilde{f}^{\lambda})) \bigcap \mathcal{U}_{\mbox{\rm loc}}^{\lambda} ;$

\item $\mathcal{T}^h_{\mbox{\rm loc}}(\mbox{\rm trop}(\widehat{f}^{\lambda})) \bigcap \mathcal{U}_{\mbox{\rm loc}}^{\lambda} ;$

\item $\overline{\mbox{\rm ord}(V( \langle \widehat{f}^{\lambda} \rangle ) )} \bigcap \mathcal{U}_{\mbox{\rm loc}}^{\lambda}. $
\end{enumerate}
\end{proposition}

Now, we have the following theorem:

\begin{theorem}\label{2-2-22} The following subsets of $\mathcal{U}_{\mbox{\rm loc}}$ coincide:

\begin{enumerate}
\item $\{ w \in \mathcal{U}_{\mbox{\rm loc}} \ | \ \forall f \in I, \ \mbox{\rm in}_w(f) \in  \mbox{\rm Gr}^w(\Ohat)  \mbox{ is not a monomial} \};$

\item $\mathcal{T}^h_{\mbox{\rm loc}}(\mbox{\rm trop}(f));$

\item $ \bigcup_{\lambda \in \Lambda} \Bigl( \mathcal{T}^h_{\mbox{\rm loc}}(\mbox{\rm trop}(\widetilde{f}^{\lambda})) \bigcap \mathcal{U}_{\mbox{\rm loc}}^{\lambda} \Bigr);$

\item $\bigcup_{\lambda \in \Lambda} \Bigl( \mathcal{T}^h_{\mbox{\rm loc}}(\mbox{\rm trop}(\widehat{f}^{\lambda})) \bigcap \mathcal{U}_{\mbox{\rm loc}}^{\lambda} \Bigr);$

\item $\bigcup_{\lambda \in \Lambda} \Bigl( \overline{\mbox{\rm ord}(V( \langle \widehat{f}^{\lambda} \rangle ) )} \bigcap \mathcal{U}_{\mbox{\rm loc}}^{\lambda} \Bigr).$
\end{enumerate}

\end{theorem}

In view of this theorem, we will define our local tropical variety
of $f$ as the union of the tropical varieties defined on strata 
$\mathcal{U}_{\rm loc}^\lambda$
for finite polynomials $\hat f$.

\begin{definition}\label{2-2-3}  Let $I= \langle f \rangle_{\Ohat}$ be a principal ideal for some $f \in \Ohat$. 
Then we define the local tropical variety of $I$ as 
$$ \mathcal{T}_{\mbox{\rm loc}}(I ) : = \bigcup_{\lambda \in \Lambda} \Bigl( \overline{\mbox{\rm ord}(V( \langle \widehat{f}^{\lambda} \rangle ) )} \bigcap \mathcal{U}_{\mbox{\rm loc}}^{\lambda} \Bigr).$$
\end{definition}

The following proposition shows that our definition is compatible with
that of polynomial tropical varieties.

\begin{proposition}\label{2-2-4} For $f \in \mathbb{C}[x] \subset \Ohat$, let $I= \langle f \rangle_{\mathbb{C}[x]} \mbox{ and } I^e= \langle f \rangle_{\Ohat}$.
Then $\mathcal{T}(I) \bigcap \mathcal{U}_{\mbox{\rm loc}}$ 
and $\mathcal{T}_{\mbox{\rm loc}}(I^e)$ are equal as sets.
\end{proposition}

\begin{proof}
It follows from Theorem \ref{main} that we have
$\mathcal{T}(I) \bigcap \mathcal{U}_{\mbox{loc}} = \mathcal{T}_h(\mbox{trop}(f)) \bigcap \mathcal{U}_{\mbox{loc}} =  \mathcal{T}^h_{\mbox{loc}}(\mbox{trop}(f))$.
By Theorem \ref{2-2-22}, 
$$ \mathcal{T}^h_{\mbox{loc}}(\mbox{trop}(f))=\bigcup_{\lambda \in \Lambda} \Bigl( \overline{\mbox{ord}(V( \langle \widehat{f}^{\lambda} \rangle ) )} \bigcap \mathcal{U}_{\mbox{loc}}^{\lambda} \Bigr) = \mathcal{T}_{\mbox{loc}}(I^e ).$$
\end{proof}

We have introduced auxiliary polynomials $\hat f$
to define a local tropical variety 
with utilizing results on polynomial tropical varieties.
It is an open question to give a description 
of our local tropical zero set without introducing the auxiliary polynomials $\hat f$.
For example, we conjecture that we do not need $\hat f$ 
to define our local tropical zero set on the maximal stratum.
In fact, 
let us introduce the set of Puiseux series of positive order 
$$ K_+ = \{ p(t) \in K=\C\{\{ t\}\} \,|\, {\rm ord}(p) > 0 \}.$$
Note that the composition of $f$ and $p \in K_+^n$ is well-defined.
We conjecture that
$$\overline{\mbox{\rm ord}(V( \langle {\tilde f}^{0} \rangle ) )}
 \bigcap \mathcal{U}_{\mbox{\rm loc}}^{0} 
= \overline{\mbox{\rm ord}(V_+(f))}\bigcap \mathcal{U}_{\mbox{\rm loc}}^{0} 
$$ 
where
$$ V_+(f) = \{ p(t)=(p_1(t), \ldots, p_n(t)) \in K_+ ^n \,|\, 
   f(p) = 0 \}.
$$

\section{Local tropical variety (general case)}

Suppose that $I \subset \Ohat$ is any ideal.
In this section, most of the proofs follow from the principal ideal case.

\subsection{Tropical Variety on the Maximal Stratum}

In this subsection, we suppose that $w$ lies in $\mathcal{U}_{\mbox{\rm loc}}^0 :=\{ (u_1,\ldots,u_n) \ | \ \forall i, u_i > 0 \}$.

For each $f\in I$, similarly to Proposition \ref{E(f)}, we define $E_0(f) \subset \mbox{Supp}(f)$. 
Then, we set $\widetilde{f}^0 := \sum_{\beta \in E_0(f) }a_{\beta}x^{\beta} \in \mathbb{C}[x]$.

\begin{proposition}\label{2-3-4} The following subsets of $\mathcal{U}_{\mbox{\rm loc}}^0$  coincide:
\begin{enumerate}
\item $\{ w \in \mathcal{U}_{\mbox{\rm loc}}^0 \ | \ \mbox{For any } f \in I, \mbox{\rm in}_w(f) \in \mbox{\rm Gr}^w(\Ohat)  \mbox{is not a monomial} \};$ 

\item $\big ( \bigcap_{f \in I }\mathcal{T}^h_{\mbox{\rm loc}}(\mbox{\rm trop}(f)) \big ) \bigcap \mathcal{U}_{\mbox{\rm loc}}^0 ;$

\item $\big ( \bigcap_{f \in I }\mathcal{T}^h_{\mbox{\rm loc}}(\mbox{\rm trop}(\widetilde{f^0})) \big ) \bigcap \mathcal{U}_{\mbox{\rm loc}}^0 ;$ 

\item $\Bigl( \bigcap_{f \in I } \overline{\mbox{\rm ord}(V(  \langle \widetilde{f^0}  \rangle  ) )} \Bigr) \bigcap \mathcal{U}_{\mbox{\rm loc}}^0 .$

\end{enumerate}
\end{proposition}

\subsection{Tropical Variety on a substratum $\mathcal{U}_{\mbox{\rm loc}}^1$}

In this subsection, we suppose that $w$ lies in $\mathcal{U}_{\mbox{\rm loc}}^1 :=\{ (0,u_2,\ldots,u_n) \ | \ \forall i ( \neq 1 ), u_i > 0 \}$.

Similarly to Definitions \ref{25--1} and \ref{2-2-16}, for each $f\in I$, we define:

\begin{definition}\label{2-2-7} $\widetilde{f}^1 := \sum_{\beta \in E_1(f) }a_{\beta}x^{\beta}$ and 
$\widehat{f}^1 := \sum_{\beta \in A_1(f) }a_{\beta}x^{\beta}.$
\end{definition}

\begin{proposition}\label{2-3-5} The following subsets of $\mathcal{U}_{\mbox{\rm loc}}^1$ coincide:
\begin{enumerate}
\item $\{ w \in \mathcal{U}_{\mbox{\rm loc}}^1 \ | \mbox{For any } f \in I, \mbox{\rm in}_w(f) \in \mbox{\rm Gr}^w(\Ohat)  \mbox{is not a monomial} \};$ 

\item $\Bigl( \bigcap_{f \in I} \mathcal{T}^h_{\mbox{\rm loc}}(\mbox{\rm trop}(f)) \Bigr) \bigcap \mathcal{U}_{\mbox{\rm loc}}^1 ;$

\item $\Bigl( \bigcap_{f \in I} \mathcal{T}^h_{\mbox{\rm loc}}(\mbox{\rm trop}(\widetilde{f^1})) \Bigr)  \bigcap \mathcal{U}_{\mbox{\rm loc}}^1 ;$

\item $\Bigl( \bigcap_{f \in I} \mathcal{T}^h_{\mbox{\rm loc}}(\mbox{\rm trop}(\widehat{f^1})) \Bigr)  \bigcap \mathcal{U}_{\mbox{\rm loc}}^1 ;$

\item $ \Bigl( \bigcap_{f\in I}\overline{\mbox{\rm ord}(V(  \widehat{f^1}  ) )} \Bigr) \bigcap \mathcal{U}_{\mbox{\rm loc}}^1 .$
\end{enumerate}
\end{proposition}

\subsection{Tropical Variety on General Strata}

In this subsection, we suppose that $w$ lies in $\mathcal{U}_{\mbox{\rm loc}}$.

We will use the same argument as in 4.2. For 
$\mathcal{U}_{\mbox{loc}}^j :=\{ (u_1,\ldots,u_n) \ | \ \forall i ( \neq j ), u_i > 0, u_j=0\}, \ldots, \mathcal{U}_{\mbox{loc}}^n, \mathcal{U}_{\mbox{loc}}^{12} :=\{ (0,0,u_3\ldots,u_n) \ | \ \forall i ( \neq 1,2 ), u_i > 0 \},\ldots,\mathcal{U}_{\mbox{loc}}^{234\cdots n} :=\{ (0,\ldots,0,u_n) \ | \  u_n > 0, \forall j ( \neq n ),u_j=0\}, \mathcal{U}_{\mbox{loc}}^{123\cdots n} :=\{ (0,\ldots,0) \}$ .

\begin{definition} For $\mathcal{U}_{\mbox{loc}}^{123\cdots n} :=\{ (0,\ldots,0) \}$,
we define $\widetilde{f}^{12\cdots n}=\widehat{f}^{12\cdots n} : = 0$.
For $\mathcal{U}_{\mbox{loc}}^{0} :=\{ (u_1,\ldots,u_n) \ | \ \forall i, u_i > 0 \}$, 
we define $\widetilde{f}^{0}=\widehat{f}^{0} $.
\end{definition}
Then, for each $\lambda \in \Lambda = \{ 0,1,2,\ldots,n,12,13,\ldots,234\cdots n,123\cdots n \}$, 
the following proposition holds:

\begin{proposition}\label{2-3-8} The following subsets of $\mathcal{U}_{\mbox{\rm loc}}^{\lambda}$ coincide:
\begin{enumerate}
\item $\{ w \in \mathcal{U}_{\mbox{\rm loc}}^{\lambda} \ | \mbox{For any } f \in I, \mbox{\rm in}_w(f) \in \mbox{\rm Gr}^w(\Ohat)  \mbox{is not a monomial} \};$ 
\item $\Bigl( \bigcap_{f \in I} \mathcal{T}^h_{\mbox{\rm loc}}(\mbox{\rm trop}(f)) \Bigr) \bigcap \mathcal{U}_{\mbox{\rm loc}}^{\lambda}; $
\item $\Bigl( \bigcap_{f \in I} \mathcal{T}^h_{\mbox{\rm loc}}(\mbox{\rm trop}(\widetilde{f^{\lambda}})) \Bigr)  \bigcap \mathcal{U}_{\mbox{\rm loc}}^{\lambda}; $
\item $\Bigl( \bigcap_{f \in I} \mathcal{T}^h_{\mbox{\rm loc}}(\mbox{\rm trop}(\widehat{f^{\lambda}})) \Bigr)  \bigcap \mathcal{U}_{\mbox{\rm loc}}^{\lambda} ;$
\item $ \Bigl( \bigcap_{f\in I}\overline{\mbox{\rm ord}(V(  \widehat{f^{\lambda}}  ) )} \Bigr) \bigcap \mathcal{U}_{\mbox{\rm loc}}^{\lambda} .$
\end{enumerate}
\end{proposition}

\begin{theorem}\label{2-2-2} The following subsets of $\mathcal{U}_{\mbox{\rm loc}}$ coincide:
\begin{enumerate}
\item $\{ w \in \mathcal{U}_{\mbox{\rm loc}} \ | \ \mbox{For any } f \in I, \mbox{\rm in}_w(f) \in \mbox{\rm Gr}^w(\Ohat)  \mbox{is not a monomial} \};$ 
\item $\bigcap_{f\in I}\mathcal{T}^h_{\mbox{\rm loc}}(\mbox{\rm trop}(f))  ;$
\item $\bigcup_{\lambda \in \Lambda} \Biggl( \Bigl( \bigcap_{f \in I} \mathcal{T}^h_{\mbox{\rm loc}}(\mbox{\rm trop}(\widetilde{f^{\lambda}})) \Bigr)  \bigcap \mathcal{U}_{\mbox{\rm loc}}^{\lambda} \Biggr) ;$
\item $\bigcup_{\lambda \in \Lambda} \Biggl( \Bigl( \bigcap_{f \in I} \mathcal{T}^h_{\mbox{\rm loc}}(\mbox{\rm trop}(\widehat{f^{\lambda}})) \Bigr)  \bigcap \mathcal{U}_{\mbox{\rm loc}}^{\lambda} \Biggr) ;$
\item $\bigcup_{\lambda \in \Lambda} \Biggl( \Bigl( \bigcap_{f\in I}\overline{\mbox{\rm ord}(V(  \widehat{f^{\lambda}}  ) )} \Bigr) \bigcap \mathcal{U}_{\mbox{\rm loc}}^{\lambda} \Biggr). $
\end{enumerate}
\end{theorem}

\begin{definition}\label{2-2-10}  Let $I \subset \Ohat$ be an ideal. 
We define the local tropical variety of $I$ as
$$ \mathcal{T}_{\mbox{\rm loc}}(I ) : = \bigcup_{\lambda \in \Lambda} \Biggl( \Bigl( \bigcap_{f\in I}\overline{\mbox{\rm ord}(V(  \widehat{f^{\lambda}}  ) )} \Bigr) \bigcap \mathcal{U}_{\mbox{\rm loc}}^{\lambda} \Biggr). $$
\end{definition}

\begin{proposition}\label{2-2-11} Let $I \subset \mathbb{C}[x] \subset \Ohat$ be an ideal 
and $ I^e= I \cdot \Ohat$.
Then $\mathcal{T}(I) \bigcap \mathcal{U}_{\mbox{\rm loc}}$ 
and $\mathcal{T}_{\mbox{\rm loc}}(I^e)$ are equal.
\end{proposition}

\begin{proof}
We show the following two statements:

For each $\lambda \in \Lambda$, 
\begin{enumerate}
\item[(i)] For $w\in \mathcal{U}_{\mbox{loc}}^{\lambda}$, if
there exists a $f \in I^e $ s.t. $\mbox{in}_w(f)$ being a monomial, then there exists a $\tilde{f} \in I$ s.t. $\mbox{in}_w(\tilde{f})$ being a monomial.

\item[(ii)] For $w\in \mathcal{U}_{\mbox{loc}}^{\lambda}$, if 
there exists a $f \in I $ s.t. $\mbox{in}_w(f)$ being a monomial, then there exists a $\tilde{f} \in I^e$ s.t. $\mbox{in}_w(\tilde{f})$ being a monomial.
\end{enumerate}

The proof of the proposition follows from the statements above. The second statement (ii) is clear.

Let $f=\sum_{i=1}^l h_ig_i\in I^e$, where $h_i \in \Ohat \mbox{ and } g_i \in I$. 
Suppose $\mbox{in}_w(f)=m$ is a monomial. 
Now, for each $h_i$, we take $\widehat{h_i^{\lambda}} \in \mathbb{C}[x]$ as same as in Definition \ref{2-2-16}.
This implies that $f':=\sum_{i=1}^l \widehat{h_i^{\lambda}} g_i$ is an element of $I$.  

With the notations above, the following lemma holds: 
(This Lemma proves statement (i) and completes the proof of proposition.)
\end{proof}
\begin{lemma}\label{3-6}
$\mathcal{T}_h(\mbox{\rm trop}(f') ) \bigcap \mathcal{U}_{\mbox{loc}}^{\lambda} = \mathcal{T}^h_{\mbox{loc}}(\mbox{\rm trop}(f) ) \bigcap \mathcal{U}_{\mbox{loc}}^{\lambda}.$ 
\end{lemma}
\begin{proof}
Lemma \ref{25-3} implies that, for each $w \in \mathcal{U}_{\mbox{loc}}^{\lambda}$, if $h_{ij}$ is a term of $h_i$ 
which have the smallest $w$-weight among the terms of $h_i$, then it is also an element of $\widehat{h_i^{\lambda}}$. 
That means, a term of $\mbox{\rm in}_w(f)$ is also a term of $f':=\sum_{i=1}^l \widehat{h_i^{\lambda}} g_i$.
We conclude as in the proof of Proposition \ref{25-1}.
\end{proof}

\section{Local Gr\"obner fan}

First, we will introduce the local Gr\"obner fan
following Assi, Castro, and Granger \cite{acg01}.
Let $u\in \mathcal{U}_{\mbox{loc}}$ and define
$$S(u) := \{u' \in \mathcal{U}_{\mbox{loc}}~|~\mbox{\rm Gr}^u(\Ohat)= \mbox{\rm Gr}^{u'}(\Ohat) \}.$$
Then, for a given ideal $I$ in $\Ohat$, we define the equivalance relation
$$u \sim u' \Leftrightarrow u' \in S(u) \mbox{ and } \mbox{\rm in}_u(I)= \mbox{\rm in}_{u'}(I).$$
For a local weight vector $u \in \mathcal{U}_{\mbox{loc}}$, 
we call $\mbox{supp}(u)=\{i ~|~ u_i \neq 0 \}$ the support of $u$. 
By Lemma \ref{3-13-1}, we have 
$\mbox{\rm Gr}^u(\Ohat)= \mbox{\rm Gr}^{u'}(\Ohat) \Leftrightarrow \mbox{supp}(u)=\mbox{supp}(u')$.
So we have:
$$u \sim u' \Leftrightarrow \mbox{supp}(u)=\mbox{supp}(u') \mbox{ and } \mbox{\rm in}_u(I)= \mbox{\rm in}_{u'}(I).$$ 

\begin{definition}\label{2-7-2}  We call the equivalence class:
$$\mathcal{C}[u]: = \{u' \in  \mathcal{U}_{\mbox{loc}} ~|~  u \sim u' \}$$
a local open Gr\"obner cone (the local Gr\"obner cone of $I$ w.r.t. $u$).
And we define the set of closures of equivalence classes:
$$\mbox{\rm LGF}(I):= \{ \overline{\mathcal{C}[u]} ~|~ u \in \mathcal{U}_{\mbox{loc}} \}$$
as the local Gr\"obner fan of $I$.
\end{definition}

\begin{theorem}{\rm \cite{Bah}}\label{2-7-3}
Let $I$ be an ideal in $\Ohat$. The local Gr\"obner fan $\mbox{\rm LGF}(I)$ 
is a polyhedral fan.
\end{theorem}

In this section, we will consider the relation between the local Gr\"obner fan and the local tropical variety of 
an ideal in $\Ohat$. First, we will state some Propositions.

\begin{proposition}\label{211916}    Let $I$ be an ideal in $\Ohat$.
Then the following subsets of $\mathcal{U}_{\mbox{\rm loc}}$ coincide.

\begin{enumerate}
\item $\{ w\in \mathcal{U}_{\mbox{\rm loc}} \ | \  \mbox{for any } f \in I, \mbox{\rm in}_w(f)  \in  \mbox{\rm Gr}^w(\Ohat) \mbox{ is not a monomial.} \}$

\item $\{ w\in \mathcal{U}_{\mbox{\rm loc}} \ | \ \mbox{\rm in}_w(I) \subset \mbox{\rm Gr}^w(\Ohat)\mbox{ contains no monomial.}   \}$
\end{enumerate}
\end{proposition}
The Proposition can be proved with showing
that if ${\rm in}_w(I)$ contains a monomial, then
there exists $f \in I$ such that ${\rm in}_w(f)$ is a monomial.

Let $I$ be an ideal in and $\mathcal{C}[w]$ be a local open Gr\"obner cone of $I$.
For a local vector $w' \in \overline{\mathcal{C}[w]} \setminus \mathcal{C}[w]$, 
put ${\tilde w}=w'+\epsilon \cdot w$ 
for some $\epsilon > 0$ sufficiently small.

The following proposition seems to be well-known in the ring of polynomials.
For the case of power series, it is called Assi's twin lemma \cite{acg01}.
Since we do not find a proof of this fact in literatures, 
we will also include a proof.
\begin{proposition}\label{1-18}  
$\mbox{\rm in}_{\tilde w}(I) = \mbox{\rm in}_w(\mbox{\rm in}_{w'}(I))$ in $\mbox{\rm Gr}^w(\Ohat)$.   
\end{proposition}
\begin{proof}
First, we will show that both
$\mbox{\rm in}_{\tilde w}(I)$ and 
$\mbox{\rm in}_w(\mbox{\rm in}_{w'}(I))$ 
are ideals of the same graded ring $\mbox{\rm Gr}^{\tilde w}(\Ohat)$ 
(i.e. $\mbox{Gr}^w(\mbox{Gr}^{w'}(\Ohat))=\mbox{Gr}^{\tilde w}(\Ohat)$).

For $\epsilon > 0$, the point $\frac{1}{1+\epsilon} (w' + \epsilon w)$ 
lies on the open segment $(w,w')$.
Then, we have $w' + \epsilon w \in C[w]$ by the convexity
and properties of cones.
(Note that \cite{Bah} proved the polyhedral property of local Gr\"obner fan
without using this proposition.)
This implies $C[w] = C[\tilde w]$ and consequently 
$\mbox{supp}(w)=\mbox{supp}(\tilde w)$.

Next, for the proof, we will show the following definition and lemmas:

\begin{definition}\label{2-7-11}\rm{[\cite{Bah} p7,8, Definition 2.1.3.]}    
Let $f$ be an element of $\Ohat$ and $\prec $ be a monomial order.
We denote by $\mbox{\rm exp}_{\prec}(f)$ the maximal element of $\mbox{\rm Supp}(f)$ w.r.t. $\prec $
 and call it leading exponent.
\begin{enumerate}
\item  Let $I$ be an ideal in $\Ohat$. We define the set of the leading exponents of $I$ as
$$\mbox{\rm Exp}_{\prec}(I)=\{ \mbox{\rm exp}_{\prec}(f) ~|~ f \in I, ~f \neq 0 \}.$$
Then, there exists $\mathcal{G}= \{ g_1,\ldots,g_r \} \subset I$ such that 
$\mbox{\rm Exp}_{\prec}(I) = \bigcup_j( \mbox{\rm exp}_{\prec}(g_j)+\mathbb{N}^n)$.
Such a set $\mathcal{G}$ is called a $\prec$-standard basis of $I$.

\item  Let $J$ be a $u$-homogeneous ideal (i.e. generated by homogeneous elements) in $\mbox{\rm Gr}^u(\Ohat)$. 
We define the set of the leading exponents of $J$ as
$$\mbox{\rm Exp}_{\prec}(J)=\{ \mbox{\rm exp}_{\prec}(f) ~|~ f \in J, ~f \neq 0 \}.$$
Then, there exists $\mathcal{G}= \{ g_1,\ldots,g_r \} \subset J$ made of homogeneous elements 
such that 
$\mbox{\rm Exp}_{\prec}(J) = \bigcup_j( \mbox{\rm exp}_{\prec}(g_j)+\mathbb{N}^n)$.
Such a set $\mathcal{G}$ is called a (homogeneous) $\prec$-standard basis of $J$.
\end{enumerate}

Now given a standard basis $G$ of $J \subset \mbox{\rm Gr}^u(\Ohat)$ where $u \in \mathcal{U}_{\mbox{\rm loc}}$.
 We say that $G$
is minimal if for $g,g' \in G$, $\exp_\prec(g) \in \exp_\prec(g')+\N^n$
implies $g=g'$. We say that $G$ is reduced if it is minimal,
unitary (i.e. $\lc_\prec(g)=1$ for any $g\in G$) and if for any
$g\in G$, $\mbox{\rm Supp}(g) \smallsetminus \{\mbox{\rm exp}_\prec(g)\} \subset \N^n
\smallsetminus \mbox{\rm Exp}_\prec(J)$.
\end{definition}

\begin{lemma}\label{21186}\rm{[\cite{Bah} Lemma 2.1.5.]}  
Given $I$ in $\Ohat$ and $u \in \mathcal{U}_{\mbox{\rm loc}}$, if $G$ is 
the reduced $\prec_u$-standard basis
then $\ini_u(G)$ is the reduced $\prec$-standard basis of $\ini_u(I)$.
\end{lemma}

\begin{lemma}\label{211811}\rm{[\cite{Bah} Lemma 2.3.1.]} 
Given an ideal $I$ in $\Ohat$ and $u \in \mathcal{U}_{\mbox{\rm loc}}$, 
for any $u'\in \overline{C_I[u]} \smallsetminus C_I[u]$ there exists a
local order $\prec$ such that the reduced standard bases of $I$ with
respect to $\prec_u$ and to $\prec_{u'}$ are the same. 
(In this proof, we have $\prec = \prec^1_u$, where $\prec^1 \mbox{ is a local order}$  satisfying the statements.)
\end{lemma}

With these lemmas, we will prove Proposition \ref{1-18}.

[Proof of Proposition \ref{1-18} (continued)] 
First, we take $\prec=\prec^1_w$,  $\prec^1$ being a local order. 
Clearly, $\prec_w=\prec$.
Let $\mathcal{G}$ be the $\prec_u$-reduced standard basis of $I$.
Then, by Lemma \ref{21186}, $\mbox{in}_w(\mathcal{G})$ is the reduced $\prec$-standard basis of $\mbox{\rm in}_w(I)$.
And, by Lemma \ref{211811}, $\mathcal{G}$ is also the reduced $\prec_{w'}$-standard basis of $I$.
By Lemma \ref{21186}, it implies that $\mbox{in}_{w'}(\mathcal{G})$ is the reduced $\prec$-standard basis of $\mbox{in}_{w'}(I)$.
Since $\prec_w=\prec$, $\mbox{in}_w(\mbox{in}_{w'}(\mathcal{G}))$ is the reduced $\prec$-standard basis of 
$\mbox{in}_w(\mbox{in}_{w'}(I))$, by Lemma \ref{21186}.
Since $w=w'+\epsilon \cdot w$, 
we have $\mbox{in}_w(\mbox{in}_{w'}(\mathcal{G})) =\mbox{in}_{w}(\mathcal{G})$.
Thus, both $\mbox{in}_w(I)$ and $\mbox{in}_w(\mbox{in}_{w'}(I))$ are ideals of 
the same ring and have the same basis. This proves the statement.
\end{proof}

Now we have the following theorem:
\begin{theorem}\label{211852} Let $I$ be an ideal in $\Ohat$.
Then the following subset of $\mathcal{U}_{\mbox{\rm loc}}$, 
the local tropical variety of $I$,
is a subfan of
the local Gr\"obner fan of $I$:
$$\mathcal{T}_{\mbox{\rm loc}}(I) = \{ w \in \mathcal{U}_{\mbox{\rm loc}} \ | \ \mbox{\rm in}_w(I) \subset  \mbox{\rm Gr}^w(\Ohat) \mbox{ contains no monomials} \}.$$
\end{theorem}
\begin{proof}
Let  
$\mbox{\rm LGF}(I) = \{ \overline{\mathcal{C}_1},\ldots,\overline{\mathcal{C}_r} \}$ 
where $\mathcal{C}_1, \ldots, \mathcal{C}_r$ are the open Gr\"obner cones.
Suppose that, for any vector $w$ in  $\mathcal{C}_i$ ($1 \leq i \leq l$), 
$\mbox{in}_w(I)$ contains no monomial, and, 
for any vector $w$ in  $\mathcal{C}_j$ ($l+1 \leq j \leq r$), 
$\mbox{in}_w(I)$ contains a monomial.
Then for any vector $\overline{C_i} \setminus C_i \ni w'$ ($1 \leq  i \leq l$), 
$\mbox{in}_{w'}(I)$ contains no monomial
(proof: $C_i \ni w = w'+\epsilon a$ for some 
 $a \in \mathbb{R}^n $ and $ \epsilon > 0$ is sufficiently small.
Then we have $\mbox{in}_w(I)=\mbox{in}_a(\mbox{in}_{w'}(I))$.
If $\mbox{in}_{w'}(I)$ contains a  monomial, then $\mbox{in}_{w}(I)$ also contains a monomial. 
This contradicts to hypothesis.).
Thus, we have $\mathcal{T}_{\mbox{\rm loc}}(I) \supset \{ \overline{C_1},\ldots,\overline{C_l} \}$.
Now we suppose that $\mbox{in}_u(I)$ contains no monomial
for some $u \in \overline{C_j} \setminus C_j ~(l+1 \leq j \leq r)$.
Then, by finiteness of the number of possible inital ideals, we have $u \in C_i \mbox{ for some } 1\leq i \leq l$.
Thus we have $\mathcal{T}_{\mbox{\rm loc}}(I) = \{ \overline{C_1},\ldots,\overline{C_l} \}$.
Since $\mbox{\rm LGF}(I)$ is a polyhedral fan, 
the set $\{ \overline{C_1},\ldots,\overline{C_l} \}$ also satisfies the properties of a polyhedral fan.
\end{proof}

\section{Tropical finite set}

In this section, we want to find a finite subset $\mathcal{H}$ of an ideal $I$ in $\Ohat$ satisfying
$$\bigcap_{f\in I}\mathcal{T}^h_{\mbox{loc}}(\mbox{trop}(f)) = \bigcap_{f\in \mathcal{H}}\mathcal{T}^h_{\mbox{loc}}(\mbox{trop}(f)).$$
For this purpose, we will construct a finite subset $\mathcal{H}$ of  $I$ satisfying the following condition:

\begin{equation}\mbox{"For each } w \in \mathcal{U}_{\mbox{loc}}, \mbox{ if } w \notin \bigcap_{f\in I}\mathcal{T}^h_{\mbox{loc}}(\mbox{trop}(f)), \tag{$\ast$}
\end{equation} 
$$\mbox{ then } \{\mbox{in}_w(h) \in \mbox{Gr}^w(\Ohat) \ | \ h \in \mathcal{H} \} \mbox{ contains a monomial."}$$

\

To prove the existence of such a set, let us state lemmas and propositions.

\begin{lemma}\label{2-1-1} Let $h$ be in $\Ohat$. 
Then $h$ is $w$-homogeneous if and only if 
$ \{ w \cdot \alpha \ | \ \alpha \in \mbox{\rm Supp}(h) \} = \{c \}$ for some $c \in \mathbb{R}$. 
\end{lemma}

\begin{lemma}\label{2-1-2} Given an ideal $I$ in 
$\Ohat$, let $\mathcal{C}[w]$ be a local open Gr\"obner cone of $I$. Suppose 
the dimension of $\mathcal{C}[w] = k \leq n$.
Then, for any $h \in \mbox{\rm Gr}^w(\Ohat)$, the following statements are equivalent:
\begin{enumerate}
\item  $ h$ is $w'$-homogeneous, $\forall w' \in \mathcal{C}[w];$

\item  $ h$ is $(\widetilde{w}^1,\ldots, \widetilde{w}^k)$-homogeneous, where $\widetilde{w}^1, \ldots, \widetilde{w}^k$ form the  $1-\mbox{skeleton}$ of $\mathcal{C}[w];$ 

\item $ h$ is $(w^1, \ldots,  w^k)$-homogeneous, where we suppose $w^1, \ldots,  w^k\in \mathcal{C}[w]$ are independent over $\mathbb{R}.$
\end{enumerate}
\end{lemma}
\begin{proof}
It is clear that (1) implies (2).

(3) $\Rightarrow$ (1): By the hypothesis and Lemma \ref{2-1-1}, there exists 
$c_i \in \mathbb{R}$ s.t. 
$\{ w^i \cdot \alpha \ | \ \alpha \in \mbox{Supp}(h) \} = \{c_i \}$ for each $1\leq i \leq k$.
Fix $w' \in \mathcal{C}[w]$. There exist  $a_i \in \mathbb{R}$ 
s.t. $w'=\sum_{i=1}^k a_i w^i$. For any $\alpha \in \mbox{Supp}(h)$, we have 
$\alpha \cdot w' = \sum_{i=1}^k a_i (w^i \cdot \alpha)=\sum_{i=1}^k a_i c_i \in \mathbb{R}$.
By Lemma \ref{2-1-1} again, we have $ h$ is $w'$-homogeneous.

(2) $\Rightarrow$ (1): By the hypothesis, for $w' \in \mathcal{C}[w]$ there exist  $a_i \in \mathbb{R}_{> 0}$ 
s.t. $w'=\sum_{i=1}^k a_i \widetilde{w}^i$. Then the proof is similar to above proof.

(1) $\Rightarrow$ (2): Now we will prove that (1) implies the following statement (which is more general than (2)):
$\mbox{ for any } \widetilde{w} \in \overline{\mathcal{C}[w]} (\setminus \mathcal{C}[w]), h \mbox{ is } \widetilde{w}\mbox{-homogeneous}$.
Let 
$\{ w(l)\in \mathcal{C}[w] \ | \ l \in \mathbb{N} \}$ be a sequence of $\mathcal{C}[w]$ s.t. 
$\lim_{l \to \infty} w(l)  = \widetilde{w}$.
Lemma \ref{2-1-1} implies that, for any $l \in \mathbb{N}$ 
we have $w(l) \cdot (\alpha - \alpha') = 0$ where $\alpha, \alpha' \in \mbox{Supp}(h)$ are arbitrary. 
The continuity of $x \mapsto x \cdot (\alpha - \alpha')$ implies 
$\widetilde{w} \cdot (\alpha - \alpha') = 0$. 
Then $\{ \widetilde{w} \cdot \alpha \ | \ \alpha \in \mbox{Supp}(h) \} = \{ c \}$ for some $c \in \mathbb{R}$
 and, again, by Lemma \ref{2-1-1}, $h$ is $\widetilde{w}$-homogeneous.
\end{proof}

\begin{proposition}\label{2-1-3}  Let $I$ be an ideal in 
$\Ohat$ and  $\mathcal{C}[w]$ be a local open Gr\"obner cone of $I$. 
Suppose that, for $w \in \mathcal{C}[w]$, 
$\mbox{\rm in}_w(I) \subset \mbox{\rm Gr}^w(\Ohat)$ contains a monomial $m$.
Then there exist $f \in I$ s.t. $\mbox{\rm in}_{w'}(f)=m$ 
 for any $w' \in \mathcal{C}[w]$.
\end{proposition}

Let the local Gr\"obner fan of $I$ be 
$$\mbox{ \rm LGF}(I) = \{ \overline{\mathcal{C}_1[w_1]},\ldots,\overline{\mathcal{C}_r[w_r]} \}  \subset \mathcal{U}_{\mbox{\rm loc}},$$ 
where $\mathcal{C}_1[w_1], \ldots, \mathcal{C}_r[w_r]  \mbox{ are the open Gr\"obner cones}$.
Suppose that,
for any $1 \leq i \leq l$, $\mbox{in}_{w_i}(I) \subset \mbox{Gr}^{w_i}(\Ohat)$
contains no monomial, and, 
for each $ l+1 \leq j \leq r$, $\mbox{in}_{w_j}(I) \subset \mbox{Gr}^{w_j}(\Ohat)$
contains a monomial $m_j$.  
Then, by Proposition\ref{2-1-3},
for each $m_j$ we can find $f_j \in I$ s.t. $\mbox{in}_{w_j'}(f_j)=m_j$ for any 
$w_j' \in \mathcal{C}[w_j]$.
We define $\mathcal{H}= \{ f_j \ | \ l+1\leq j \leq r  \}$.

\begin{proposition}\label{2-1-4}
$\mathcal{H}= \{ f_j \ | \ l+1\leq j \leq r  \}$ is a local tropical finite set of $I$, i.e. satisfies $(\ast)$.
\end{proposition}

\begin{proof}
By Theorem \ref{2-2-2}, Proposition \ref{211916} and Proposition \ref{211852}, 
the following subsets of $\mathcal{U}_{\mbox{loc}}$ coincide:
\begin{enumerate}

\item $(\bigcap_{f\in I} \mathcal{T}^h_{\mbox{loc}}(\mbox{trop}(f)));  $

\item $\{ w \in \mathcal{U}_{\mbox{loc}} \ | \ \mbox{for any } f \in I, \ \mbox{\rm in}_w(f)  \in \mbox{Gr}^{w}(\Ohat) \mbox{ is not a monomial}  \};$

\item $\{ w \in \mathcal{U}_{\mbox{loc}} \ | \ \mbox{in}_w(I) \subset \mbox{Gr}^{w}(\Ohat)\mbox{ contains no monomial }  \};$ 

\item $ \{ \overline{\mathcal{C}_1[w_1]},\ldots,\overline{\mathcal{C}_l[w_l]} \} .$

\end{enumerate}

Let $w \notin \bigcap_{f\in I} \mathcal{T}^h_{\mbox{loc}}(\mbox{trop}(f))$,
then  $w \in \{ \mathcal{C}_{l+1}[w_{l+1}],\ldots,\mathcal{C}_r[w_r] \}$ and  the set $\{ \mbox{in}_w(h) \ | \ h \in \mathcal{H} \}$ contains a monomial.
\end{proof}
Now we have following theorem:

\begin{theorem}\label{2-1-6} Given an ideal $I$ in $\Ohat$, the following subsets of 
$\mathcal{U}_{\mbox{\rm loc}}$  coincide: 

\begin{enumerate}
\item  $\{ w \in \mathcal{U}_{\mbox{\rm loc}} \ | \ \mbox{ for any } f \in I, \ \mbox{\rm in}_w(f)  \in \mbox{\rm Gr}^{w}(\Ohat) \mbox{ is not a  monomial} \};$

\item $\{ w \in \mathcal{U}_{\mbox{\rm loc}} \ | \ \mbox{\rm in}_w(I) \subset \mbox{\rm Gr}^{w}(\Ohat)\mbox{ contains no monomial}  \};$ 
 
\item $ \bigcap_{f \in I }\mathcal{T}^h_{\mbox{\rm loc}}(\mbox{\rm trop}(f)) ;$

\item $\mathcal{T}_{\mbox{\rm loc}}(I);$

\item $ \{ \overline{\mathcal{C}_1[w_1]},\ldots,\overline{\mathcal{C}_l[w_l]} \} ;$

\item $ \bigcap_{l+1\leq j \leq r }\mathcal{T}^h_{\mbox{\rm loc}}(\mbox{\rm trop}(f_j)). $
\end{enumerate}
\end{theorem}

Finally, we will finish this paper with an example.

\begin{example} \rm
Let $f_1=1+x+\frac{x^2}{2}+\frac{x^3}{6}+\cdots$ be the Maclaurin's expansion of $e^x$.
And let
$f=y(f_1-1)-x^2=y(1+x+\frac{x^2}{2}+\frac{x^3}{6}+\cdots -1)-x^2=xy-x^2+\frac{x^2y}{2}+\frac{x^3y}{6}+\cdots$.
Now we will compute the local tropical variety of $\langle f \rangle$.

First, as in  Definitions \ref{25-0} and  \ref{2-2-16}, 
for $\mathcal{U}_{\mbox{ \rm loc}}^0=\{ (w_1,w_2) ~|~ w_1,w_2 > 0 \}$ we have $\widetilde{{f}^0}=\widehat{{f}^0}= xy-x^2$,
for $\mathcal{U}_{\mbox{ \rm loc}}^1=\{ (0,w_2) ~|~ w_2 > 0 \}$ we have $\widehat{{f}^1}=-x^2$,
for $\mathcal{U}_{\mbox{ \rm loc}}^2=\{ (w_1,0) ~|~ w_1 > 0 \}$ we have $\widehat{{f}^2}=xy$,
and for $\mathcal{U}_{\mbox{ \rm loc}}^{12}=(0,0)$ we have $\widehat{{f}^{12}}=0$.
By Definition \ref{2-2-3}, the local tropical variety of $\langle f \rangle$ is as in 
figure \rm{\ref{ltv1}}.

\begin{figure}[ht]
		\begin{center}
		\reflectbox{\includegraphics[width=6cm,clip]{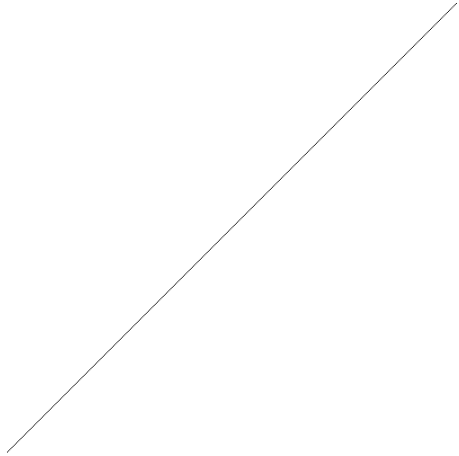}}
		\end{center}
		\caption{the local tropical variety of $f$ }
		\label{ltv1}
	\end{figure}

\end{example}


\newpage


\begin{thebibliography}{xxxxxxxxxx}


\bibitem{America}
The American Institute of Mathematics: 
Amoebas and tropical geometry, 2004.

\bibitem{acg01}
A.~Assi, F.J.~Castro-Jim\'enez, M.~Granger,
The analytic standard fan of a $\mathcal{D}$-module,
Journal of Pure and Applied Algebra {\bf 164} (2001), 3--21.


\bibitem{Bah}
R.~Bahloul and N.~Takayama:
Local Gr\"obner fan: polyhedral and computational approach.
arXiv:math.AG/0412044, 2004.


\bibitem{CLO}
D.~Cox, J.~Little, and D.~O'Shea:
Using algebraic geometry.
Number 185 in Graduate Texts in Mathematics. Springer-Verlag, New York, 1998.

\bibitem{Eis}
D.~Eisenbud: Commutative algebra with a view toward algebraic geometry. Graduate Texts in Mathematics, 150. Springer-Verlag, New York, 1995.

\bibitem{Kap.1}
M.~Einsiedler, M.~Kapranov, and D.~Lind: Non-archimedean amoebas and tropical varieties. arXiv:math.AG/0408311, 2004.


\bibitem{GKZ}
I.~Gelfand, M.~Kapranov and A.~Zelevinsky:
Discriminants, resultants and multidimensional determinants.
Birkh\"{a}user, Boston, 1994.

\bibitem{Sing}
G.~Greuel and G.~Pfister: A Singular introduction to commutative algebra. 
Springer-Verlag, Berlin, 2002.


\bibitem{Mik.1}
G.~Mikhalkin: Amoebas of algebraic varieties and tropical geometry. 
arXiv:math.AG/0403015, 2004.

\bibitem{Pen.1}
H.~Pennaneac'h: Tropical geometry and amoebas, 2003.

\bibitem{Stur.0}
J.~Richter-Gebert, B.~Sturmfels and T.~Theobald: First steps in tropical geometry. arXiv:math.AG/0306366, 2003.


\bibitem{Spey}
D.~Speyer and B.~Sturmfels: The tropical Grassmannian, Advances in Geometry,
{\bf 4} (2004), 389--411. 
arXiv:math.AG/0304218.




\bibitem{Stur.1}
B.~Sturmfels: Gr\"obner bases and convex polytopes, University Lecture Series 8, American Mathematical Society, 1996.

\bibitem{Stur.2}
B.~Sturmfels: Solving systems of polynomial equations, CBMS Regional Conference Series in Math, 97, American Mathmatical Society, 2002.





\end{thebibliography}
\end{document}